\documentclass[11pt]{amsart}
\usepackage{amssymb}
\usepackage{amsmath}
\usepackage{amsthm}
\usepackage{latexsym}
\usepackage{amsfonts}
\usepackage{graphicx}
\usepackage{graphics}
\usepackage[T1]{pbsi}

\theoremstyle{definition}
\newtheorem*{Proof}{Proof}

\newcommand{\dis}{\displaystyle}
\newcommand{\bbb}[1]{\mbox{\boldmath$#1$}}

\newcommand{\ra}{\;\rightarrow\;}

\newcommand{\de}{\delta }

\newcommand{\e}{\varepsilon }

\newcommand{\la}{\lambda }

\newcommand{\R}{\mathbb{R}}

\newcommand{\ld}{\ldots}

\newcommand{\sm}{\smallsetminus}

\newcommand{\hs}{\hfill$\square$}

\begin{document}

\title[$A_1$ weights on $\R$, an alternative approach]{$\bbb{A_1}$ weights on $\bbb{\R}$, an alternative approach}
\author{Eleftherios N. Nikolidakis}
\footnotetext{\hspace{-0.5cm}2010 {\em Mathematics Subject Classification.} Primary 42B25; Secondary 42B99.} \footnotetext{\hspace{-0.5cm}{\em Keywords and phrases.} Muckenhoupt weights, decreasing rearrangement.}

\date{}
\maketitle
\noindent
{\bf Abstract.} We will prove that if $\phi$ belongs to the class $A_1(\R)$ with constant $c\ge1$ then the decreasing rearrangement of $\phi$, belongs to the same class with constant not more than $c$. We also find for such $\phi$ the exact best possible range of those $p>1$ for which $\phi\in L^p$. In this way we provide alternative proofs of the results that appear in \cite{1}.
\section{Introduction}\label{sec1}
\noindent

The theory of Muckenhoupt weights has been proved to be an important tool in analysis. One of the most important facts concerning these is their self improving property. A way to express this is through the so-called reverse H\"{o}lder or more generally reverse Jensen inequalities (see \cite{2}, \cite{3} and \cite{7}).

In this paper we are concerned with such weights and more precisely for those $\phi$ that belong to the class $A_1(J)$ where $J$ is an interval on $\R$. This is defined as follows

A function $\phi:J\ra\R^+$ belongs to $A_1(J)$ if there exists a constant $c\ge1$ such that the following condition is satisfied:
\begin{eqnarray}
\frac{1}{|I|}\int_I\phi(x)dx\le c\cdot ess\inf_I(\phi).  \label{eq1.1}
\end{eqnarray}
for every I subinterval of $J$, where $|\cdot|$ is the Lesbesgue measure on $\R$.

Moreover if the constant $c$ is the least for which (\ref{eq1.1}) is satisfied for any $I\subseteq J$ we say that the $A_1$ constant of $\phi$ is $c$ and is denoted by $[\phi]_1$. We say then that $\phi$ belongs to the $A_1$ class of $J$ with constant $c$ and we write $\phi\in A_1(J,c)$.

It is a known fact that if $\phi\in A_1(J,c)$ then there exists $p(c)>1$ such that $\phi\in L^p$ for every $p\in[1,p(c))$.

Moreover $\phi$ satisfies a reverse H\"{o}lder inequality for every $p\in[1,p(c))$. That is for any such $p$ there exists $C_p>1$ such that
\begin{eqnarray}
\frac{1}{|I|}\int_I\phi^p(x)dx\le C_p\bigg(\frac{1}{|I|}\int_I\phi(x)dx\bigg)^p,  \label{eq1.2}
\end{eqnarray}
for every $I$ subinterval of $J$ and every $\phi\in A_1(J,c)$.

The problem of the exact determination of the best possible constant $p(c)$ has been treated in \cite{1}. More precisely it is shown there the following:\vspace*{0.2cm} \\
\noindent
{\bf Theorem A:} {\em If $\phi\in A_1((0,1),c)$ and $c$ is greater than $1$, then $\phi\in L^p(0,1)$ for any $p$ such that $1\le p<\dfrac{c}{c-1}$. Moreover the following inequality is true
\begin{eqnarray}
\frac{1}{|I|}\int_I\phi^p(x)dx\le\frac{1}{c^{p-1}(c+p-pc)}\bigg(\frac{1}{|I|}
\int_I\phi(x)dx\bigg)^p \label{eq1.3}
\end{eqnarray}
for every $I$ subinterval of $(0,1)$ and for any $p$ in the range $[1,\dfrac{c}{c-1})$. Additionally, the constant that appears in the right of inequality (\ref{eq1.3}) is best possible.}\medskip

As a consequence of the above theorem we have that the best possible range for the $L^p$-integrability of any $\phi$ with $[\phi]_1=c$ is $[1,\dfrac{c}{c-1})$.

The approach for proving the above theorem as is done in \cite{1}, is by using the decreasing rearrangement of $\phi$ which is defined by the following equation
\begin{eqnarray}
\phi^\ast(t)=\sup_{e\subset(0,1)\atop|e|\ge t}\Big[\inf_{x\in e}\phi(x)\Big], \label{eq1.4}
\end{eqnarray}
for any $t\in(0,1]$.

Then $\phi^\ast$ is a function equimeasurable to $\phi$, non-increasing and left continuous.

The immediate step for proving Theorem A, as it appears in \cite{1} is the following:\vspace*{0.2cm}\\
\noindent
{\bf Theorem B.} {\em If $\phi\in A_1((0,1),c)$ then $\phi^\ast\in A_1((0,1),c')$ for some $c'$ such that $1\le c'\le c$.}\medskip

This is treated in \cite{1} initially for continuous functions $\phi$ and generalized to arbitrary $\phi$ by use of a covering lemma. Then applying several techniques the authors in \cite{1} were able to prove Theorem A firstly for non-increasing functions and secondly for general $\phi$ by use of Theorem B.

In this paper we provide alternative proofs of the Theorems A and B.
We first prove Theorem B without any use of covering lemmas. Then we provide a proof of Theorem A for non-increasing functions $\phi$. The proof gives in an immediate way the inequality (\ref{eq1.3}). At last we prove Theorem A in it's general form by using the above mentioned results.

Additionally, we need to say that the dyadic analogue of the above problem is solved in \cite{6} while in \cite{4} and \cite{5} related problems for estimates for the range of $p$ in higher dimensions have been treated.
\section{Rearrangements of $A_1$ weights on $(0,1)$}\label{sec2}
\noindent

We are now ready to state and prove the main theorem in this section.\vspace*{0.2cm} \\
{\bf Theorem 1:} {\em Let $\phi:(0,1)\ra\R^+$ which satisfies condition (\ref{eq1.1}) for any subinterval $I$ of $(0,1)$, and for a constant $c\ge1$. Then $\phi^\ast$ satisfies this condition with the same constant}.
\begin{Proof}
It is easy to see that in order to prove our result, we need to prove the following inequality:
\begin{eqnarray}
\frac{1}{t}\int^t_0\phi^\ast(u)du\le c\phi^\ast(t)  \label{eq2.5}
\end{eqnarray}
for any $t\in(0,1]$, due to the fact that $\phi^\ast$ is left continuous and non-increasing.

For any $\la>0$ we consider the set $E_\la=\{x\in(0,1):\phi(x)>\la\}$. Let now $\e>0$. Then we can find for any such $\e$ an open set $G_\e\subseteq(0,1)$ for which $G_{\e}\supseteq E_\la$ and $|G_{\e}\sm E_\la|<\e$. Then $G_\e$ can be decomposed as follows: $G_\e=\dis\bigcup^{+\infty}_{j=1}I_{j,\e}$, where $(I_{j,\e})$ is a family of non-overlapping open subintervals of $(0,1)$. If any two of these have a common endpoint we replace them by their union. We apply the above procedure to the new family of intervals and at last we reach to a family $(I'_{j,\e})_j$ of non-overlapping open intervals such that, if $G'_{\e}=\dis\bigcup^{+\infty}_{j=1}I'_{j,\e}$ we still have that $G_\e\supseteq E_\la$ and $|G_\e\sm E_\la|<\e$.

Additionally we have that for any $j$ such that $I'_{j,\e}\neq(0,1)$ there exists an endpoint of it such that if we enlarge this interval in the direction of this point, thus producing the interval $I'_{j,\e,\de}$ with $\de$ small enough, we have that $\underset{I'_{j,\e,\de}}{ess\inf}(\phi)\le\la$. This follows by our construction and the definition of $E_\la$. Suppose now that $|E_\la|<1$. Thus $I'_{j,\e,\de}\neq(0,1)$ for any $j,\e$ and $\de$. On each of these intervals we apply (\ref{eq1.1}). So we conclude that
\[
\frac{1}{|I'_{j,\e,\de}|}\int_{I'_{j,\e,\de}}\phi\le c\cdot\underset{I'_{j,\e,\de}}{ess\inf(\phi)}\le c\la,
\]
for every $\e,\de>0$ and $j=1,2,\ld\;.$

Letting $\de\ra0^+$ we reach to the inequality $\dfrac{1}{|I'_{j,\e}|}\dis\int_{I'_{j,\e}}\phi\le c\la$ for any $j=1,2,\ld$ and every $\e>0$.

Since $G_\e=\dis\bigcup^{+\infty}_{j=1}I'_{j,\e}$ is disjoint we must have that:
\[
\frac{1}{|G_\e|}\int_{G_\e}\phi\le\sup\bigg\{\frac{1}{|I_{j,\e}|}\int_{I_{j,\e}}
\phi: \ \ j=1,2,\ld\bigg\}\le c\la
\]
for every $\e>0$ and letting $\e\ra0^+$ we have as a result that
\[
\frac{1}{|E_\la|}\int_{E_\la}\phi\le c\la\le c\cdot \underset{E_\la}{ess\inf(\phi)}.
\]
By the definition of $E_\la$ we have that
\[
\frac{1}{|E_\la|}\int_{E_\la}\phi=\frac{1}{|E_\la|}\int^{|E_\la|}_0\phi^\ast
(u)du
\]
and of course
\[
\underset{E_\la}{ess\inf(\phi)}=\underset{(0,(E_\la)]}{ess\inf(\phi^\ast)}=
\phi^\ast(|E_\la|).
\]
since $\phi^\ast$ is left continuous.
As a consequence from the above we immediately see that
\[
\frac{1}{|E_\la|}\int^{|E_\la|}_0\phi^\ast(u)du\le c\phi^\ast(|E_\la|).
\]
The same inequality holds even in the case where $|E_\la|=1$, so $G_\e=(0,1)$.Then by relation (\ref{eq1.1}) that holds for the interval $(0,1)$
we conclude (\ref{eq2.5}). Thus we have proved that $\dfrac{1}{t}\dis\int^t_0\phi^\ast(u)du\le c\phi^\ast(t)$, for every $t$ of the form $t=|E_\la|$ for some $\la>0$.

Let now $t\in(0,1]$ and define
\[
\phi^\ast(t)=\la_1, \ \ t_1=\min\Big\{s\in(0,1]:\phi^\ast(s)=\la_1\Big\}\le t \]
Additionally $|E_{\la_1}|=t_1$. As a result
\begin{align*}
\frac{1}{t}\int^t_0\phi^\ast(u)du&\le\frac{t_1}{t}\bigg(\frac{1}{|E_{\la_1}|}
\int^{|E_{\la_1}|}_0\phi^\ast(u)du\bigg)+\frac{t-t_1}{t}\la_1\\
&\le\frac{t_1}{t}c\la_1+\frac{t-t_1}{t}\la_1\le c\la_1=\phi^\ast(t).
\end{align*}
where in the second inequality we have used the above results. Theorem 1 is now proved. \hs
\end{Proof}

We proceed now to the next section.
\section{$L^p$ integrability for $A_1$ weights on $(0,1)$}\label{sec3}
\noindent

We shall now prove the following:\vspace*{0.2cm} \\
\noindent
{\bf Theorem 2.} {\em Let $\phi\in A_1((0,1),c)$ where $c$ is greater than $1$. Then, for every $p\in[1,\dfrac{c}{c-1})$, $\phi\in L^p$ and satisfies the following inequality
\begin{eqnarray}
\frac{1}{|I|}\int_I\phi^p\le\frac{1}{c^{p-1}(c+p-pc)}\bigg(\frac{1}{|I|}\int_I\phi\bigg)^p, \label{eq3.6}
\end{eqnarray}
for every $I$ subinterval of $(0,1)$.
Moreover, inequality (\ref{eq3.6}) is best possible}.\medskip

We will need first a preliminary lemma which we state as \vspace*{0.2cm} \\
\noindent
{\bf Lemma 1.} {\em Let $g:(0,1]\ra\R^+$ be a non-increasing function. Then the following inequality is true for any $p>1$ and every $\de\in(0,1)$}
\begin{eqnarray}
\int^{\de}_0\bigg(\frac{1}{t}\int^t_0g\bigg)^pdt=-\frac{1}{p-1}\bigg(\int_0^\de g\bigg)^p\frac{1}{\de^{p-1}}+\frac{p}{p-1}\int^\de_0\bigg(\frac{1}{t}\int^t_0g\bigg)^{p-1}g(t)dt.\hspace*{-1cm}
\label{eq3.7}
\end{eqnarray}
\begin{Proof}
By using Fubini's theorem it is easy to see that
\begin{eqnarray}
\int^\de_0\bigg(\frac{1}{t}\int^t_0g\bigg)^pdt=\int^{+\infty}_{\la=0}p\la^{p-1}\bigg|\bigg\{t\in(0,\de]:
\frac{1}{t}\int^t_0g\ge\la\bigg\}\bigg|dt.  \label{eq3.8}
\end{eqnarray}
Let now $\dfrac{1}{\de}\dis\int^\de_0g=f_\de\ge f=\dis\int^1_0g$.

Then
\[
\begin{array}{l}
  \dfrac{1}{t}\dis\int^t_0g>f_\de, \ \ \forall\;t\in(0,\de) \ \ \text{while} \\ [2ex]
  \dfrac{1}{t}\dis\int^t_0g\le f_\de, \ \ \forall\;t\in[\de,1].
\end{array}
\]
Let now $\la$ be such that: $0<\la<f_\de$. Then for every $t\in(0,\de]$ we have that $\dfrac{1}{t}\dis\int^t_0g\ge\dfrac{1}{\de}\dis\int^\de_0g=f_\de>\la$. Thus
\[
\bigg|\bigg\{t\in(0,\de]:\frac{1}{t}\int^t_0g\ge\la\bigg\}\bigg|=|(0,\de]|=\de.
\]
Now for every $\la\ge f_\de$ there exists unique $a(\la)\in(0,\de]$ such that $\dfrac{1}{a(\la)}\dis\int^{a(\la)}_0g=\la$. It's existence is quaranteeded by the fact that $\la>f_\de$, that $g$ is non-increasing and that $g(0^+)=+\infty$ which may without loss of generality be assumed (otherwise we work for the $\la$'s on the interval $(0,\|g\|_\infty])$.

Then
\[
\bigg\{t\in(0,\de]:\frac{1}{t}\int^t_0g\ge\la\bigg\}=(0,a(\la)].
\]
Thus from (\ref{eq3.8}) we conclude that
\begin{align}
\int^\de_0\bigg(\frac{1}{t}\int^t_0g\bigg)^pdt&=\int^{f_\de}_{\la=0}p\la^{p-1}\cdot
\de\cdot d\la+\int^{+\infty}_{\la=f_\de}p\la^{p-1}a(\la)d\la \nonumber\\
&=\de(f_\de)^p+\int^{+\infty}_{\la=f_\de}p\la^{p-1}\frac{1}{\la}\bigg(\int^{a(\la)}_0
g(u)du\bigg)d\la,  \label{eq3.9}
\end{align}
by the definition of $a(\la)$.

As a consequence (\ref{eq3.9}) becomes
\begin{align*}
\int^\de_0\bigg(\frac{1}{t}\int^t_0g\bigg)^pdt&=\frac{1}{\de^{p-1}}\bigg(\int^\de_0g\bigg)^p
+\int^{+\infty}_{\la=f_\de}p\la^{p-2}\bigg(\int^{a(\la)}_0g(u)du\bigg)d\la \\
&=\frac{1}{\de^{p-1}}\bigg(\int^\de_0g\bigg)^p+\int^{+\infty}_{\la=f_\de}p\la^{p-2}
\bigg(\int_{\{u\in(0,\de]:\atop\frac{1}{u}\int^u_0g\ge\la\}}g(u)du\bigg)d\la\\
&=\frac{1}{\de^{p-1}}\bigg(\int^\de_0g\bigg)^p+\frac{p}{p-1}\int^\de_0g(t)\Big[
\la^{p-1}\Big]^{\frac{1}{t}\int^t_0g}_{\la=f_\de}dt\\
&=\frac{1}{\de^{p-1}}\bigg(\int^\de_0g\bigg)^p+\frac{p}{p-1}\bigg[\int^\de_0
\bigg(\frac{1}{t}\int^t_0g\bigg)^{p-1}g(t)-\bigg(\int^\de_0g(t)dt\bigg)f_\de^{p-1}\bigg]\\
&=-\frac{1}{p-1}\frac{1}{\de^{p-1}}\bigg(\int^\de_0g\bigg)^p+\frac{p}{p-1}\int^{\de}_0\bigg(
\frac{1}{t}\int^t_0g\bigg)^{p-1}g(t)dt,
\end{align*}
where in the third equality we have used Fubini's theorem and the fact that $\dfrac{1}{\de}\dis\int^\de_0g=f_\de$.

Lemma 1 is now proved. \hs
\end{Proof}

Before we prove Theorem 2 we will need the following:\vspace*{0.2cm} \\
\noindent
{\bf Lemma 2.} {\em Let $g:(0,1]\ra\R^+$ be non-increasing such that
\[
\frac{1}{t}\int^t_0g(u)du\le cg(t), \ \ \text{for every} \ \ t\in(0,1].
\]
Then for every $\de\in(0,1]$ we have the following inequality
\[
\frac{1}{\de}\int^\de_0g^p\le\frac{1}{c^{p-1}(c+p-pc)}\bigg(\frac{1}{\de}\int^\de_0g\bigg)^p,
\]
for every $p$ such that $1\le p<\dfrac{p}{p-1}$. Moreover the above inequality is sharp}.
\begin{Proof}
Fix $\de\in(0,1]$ and $p\in[1,\dfrac{c}{c-1})$. Then by Lemma 1
\begin{align}
&\int^\de_0\bigg(\frac{1}{t}\int^t_0g\bigg)^pdt=-\frac{1}{p-1}\bigg(\int^\de_0g\bigg)^p
\frac{1}{\de^{p-1}}+\frac{p}{p-1}\int^\de_0\bigg(\frac{1}{t}\int^t_0g\bigg)^{p-1}g(t)dt \nonumber \\
&\Rightarrow\frac{1}{\de}\int^\de_0\bigg[\bigg(\frac{1}{t}\int^t_0g\bigg)^{p-1}g(t)-
\frac{p-1}{p}\bigg(\frac{1}{t}\int^t_0g\bigg)^p\bigg]dt\le\frac{1}{p}\bigg(\frac{1}{\de}
\int^\de_0g\bigg)^p.  \label{eq3.10}
\end{align}
We now define the following function $h_y$, of the variable $x$ for any fixed constant $y>0$
\[
h_h(x)=x^{p-1}y-\frac{p-1}{p}x^p \ ,\ \text{for} \ \ x\in[y,cy].
\]
Then
\[
h'_y(x)=(p-1)x^{p-2}y-(p-1)x^{p-1}=(p-1)x^{p-2}(y-x)\le0, \ \ \forall\;x\ge y.
\]
Thus, $h_y$ is decreasing on the interval $[y,cy]$. We conclude that for any $x$ such that $y\le x\le cy$ we have $h_y(x)\ge h_y(cy)$.

Applying the above conclusion in the case where $x=\dfrac{1}{t}\dis\int^t_0g$, $y=g(t)$ (noting that $y\le x\le cy$, for any fixed $t$) we reach to the inequality:
\begin{align}
\bigg(\frac{1}{t}\int^t_0g\bigg)^{p-1}g(t)-\frac{p-1}{p}\bigg(\frac{1}{t}\int^t_0g\bigg)^p&\ge
c^{p-1}g^p(t)-\frac{p-1}{p}c^pg^p(t)\nonumber\\
&=c^{p-1}\bigg[1-\frac{p-1}{p}c\bigg]g^p(t), \ \ \forall\;t\in(0,1].  \label{eq3.11}
\end{align}
Applying (\ref{eq3.11}) in (\ref{eq3.10}) we have as a result that
\begin{align}
&c^{p-1}\bigg[1-\frac{p-1}{p}c\bigg]\frac{1}{\de}\int^\de_0g^p(t)dt\le\frac{1}{p}
\bigg(\frac{1}{\de}\int^\de_0g\bigg)^p \nonumber\\
&\Rightarrow\frac{1}{\de}\int^\de_0g^p\le\frac{1}{c^{p-1}[p+c-pc]}
\bigg(\frac{1}{\de}\int^\de_0g\bigg)^p, \label{eq3.12}
\end{align}
which is the inequality that is stated above.

Additionally (\ref{eq3.12}) is sharp as can be seen by using the function $g(t)=\dfrac{1}{c}t^{\frac{1}{c}-1}$, $t\in(0,1]$, for $c>1$, and $g=const$ for $c=1$.

Lemma 2 is now proved. \hs
\end{Proof}

We are now ready for the\vspace*{0.2cm} \\
\noindent
{\bf Proof of Theorem 2.} Let $I\subseteq(0,1)$ be an interval.

We set $\phi_I:I\ra\R^+$ by $\phi_I(x)=\phi(x)$, $x\in I$.

Then $\phi_I$ satisfies on $I$ the condition (\ref{eq1.1}) with constant $c$. That is $\phi\in A_1(I)$ with $A_1$-constant less or equal then $c$. Then by the results of Section \ref{sec2} and a dilation argument we conclude that
\[
\phi^\ast_I=g_I:(0,|I|]\ra\R^+ \ \ \text{satisfies}
\]
\[
\frac{1}{t}\int^t_0g_I\le cg_I(t), \ \ \text{for any} \ \ t\in(0,|I|].
\]
Then by Lemma 2 and considering the results of this Section we have the inequality:
\begin{eqnarray}
\frac{1}{t}\int^t_0g^p_I(u)du\le\frac{1}{c^{p-1}(c+p-pc)}\bigg(\frac{1}{t}\int^t_0g_I(u)du\bigg)^p \label{eq3.13}
\end{eqnarray}
for any $t\in(0,|I|]$.

By the fact now that $g_I=(\phi/I)^\ast$ and (\ref{eq3.13}) we see immediately: (for $t=|I|$) that
\[
\frac{1}{|I|}\int_I\phi^p\le\frac{1}{c^{p-1}(c+p-pc)}\bigg(\frac{1}{|I|}\int_I\phi\bigg)^p.
\]
At last, we mention that the result is best possible since Lemma 2 is proved to be sharp.

Theorem 2 is now proved. \hs
\vspace*{1cm}
\noindent
Nikolidakis Eleftherios\vspace*{0.1cm}\\
Post-doctoral researcher\vspace*{0.1cm}\\
National and Kapodistrian University of Athens\vspace*{0.1cm}\\
Department of Mathematics\vspace*{0.1cm}\\
Panepistimioupolis, GR 157 84\vspace*{0.1cm}\\
Athens, Greece \vspace*{0.1cm}\\
E-mail address:lefteris@math.uoc.gr


\begin{thebibliography}{99}
\bibitem{1} B. Bojarski, C. Sbordonc and I. Wik, ``The Muckenhoupt class $A_1(\R)$'', Studia math. 101, (2) (1992) 155-163.
%
\bibitem{2} R. Coifman and C. Fefferman, ``Weighted norm inequalities for maximal functions and singular integrals'', Studia Math. 51, (1974) 241-250.
%
\bibitem{3} F. W. Gehring, ``The $L^p$-integrability of the partial derivatives of a quasiconformal mapping'', Acta Math (130) (1973) 265-277.
%
\bibitem{4} J. Kinnunen, ``Sharp results on reverse H\"{o}lder inequalities'', Ann. Acad. Sci. Fenn. Ser A I Math. Diss. 95 (1994) 1-34.
%
\bibitem{5} J. Kinnunen, ``A stability result for Muckenhoupt's weights'', Publ. Mat. 42 (1998) 153-163.
%
\bibitem{6} A. D. Melas, ``A sharp $L^p$ inequality for dyadic $A_1$ weights in $\R^n$'', Bull. London Math. Soc. 37 (2005) 919-926.
%
\bibitem{7} B. Muckenhoupt, ``Weighted norm inequalities for the Hardy-Littlewood maximal function'', Trans. Amer. Math. Soc. 165 (1972) 207-226.
\end{thebibliography}
\end{document}